\documentclass[oneside,reqno]{amsart}
\usepackage{amsmath}
\usepackage{enumerate}
\newtheorem{theorem}{Theorem}[section]

\theoremstyle{definition}
\newtheorem{definition}[theorem]{Definition}

\numberwithin{equation}{section}

\newcommand{\blankbox}[5]

\begin{document}
\title{Boundedness of weighted multilinear $p$-adic Hardy operator on Herz type spaces}
\author{Amjad Hussain$^{1,*}$}
\author{Naqash Sarfraz$^1$}
%\date{February, 2013}
\footnote{$^{}$Department of Mathematics, Quaid-I-Azam University 45320, Islamabad 44000, Pakistan\\
$^*$Corresponding Author: ahabbasi123@yahoo.com (A. Hussain); naqashawesome@gmail.com (N. Sarfraz)}


\begin{abstract} The present article focuses on the bounds of weighted multilinear $p$-adic Hardy operators on the product of Herz spaces and Morrey-Herz spaces. The corresponding norm  on both the cases are also obtain.
\end{abstract}

\maketitle

\section{\textbf{Introduction}} The weighted Hardy operator was defined in \cite{CF} as:
$$H_{\psi}(f)(x)=\int_{0}^{1}f(tx)\psi(t)dt,$$
where $\psi:[0,1]\rightarrow [0,\infty)$ is a measurable function. Xiao \cite{X} proved the boundedness of $H_{\psi}$ on either $L^{p}(\mathbb{R}^{n}), 1\leq p\leq\infty$ or $BMO(\mathbb{R}^{n}).$ Moreover, he workout on corresponding norms as well. Fu et.al \cite{FLZ3} ensured that $H_{\psi}$ is bounded on central Morrey spaces and $\lambda$-central $BMO$ spaces. It is important to mention here that the corresponding norms were also obtained. For more detail about boundedness of $H_{\psi},$ see \cite{FLZ4,LZE}.
Interestingly, when $\psi\equiv1$ and $n=1,$ the weighted Hardy operator is reduced to classical Hardy operator \cite{H} which is defined as follows.
\begin{equation}\label{EN4}H f(x)=\frac{1}{x}\int_{0}^{x}f(t)dt,\quad   x>0,\end{equation}
 which satisfied the below inequality:
\begin{equation}\label{EN5}\|H f\|_{L^{p}(\mathbb{R}^{+})}\leq\frac{p}{p-1}\|f\|_{L^{p}(\mathbb{R}^{+})},\quad 1<p<\infty.\end{equation}
The constant $p/(p-1)$ in (\ref{EN5}) was shown sharp. The Hardy operator in higher dimensional Euclidian space $\mathbb R^n$ was given by Faris see \cite{F}, which is stated as:
\begin{equation}\label{EN6}Hf(\mathbf{x})=\frac{1}{|\mathbf{x}|^{n}}\int_{|\mathbf{y}|\leq|\mathbf{x}|}f(\mathbf{y})d\mathbf{y}.\end{equation}
Over the years Hardy operator has gained significant amount of attention due to its boundedness properties. For complete understanding of Hardy type operators we refer the reader to see \cite{CG,FLLW,FLZ,GZ,X} and the references therein.

In the past few decades, $p$-adic analysis gained impeccable attraction in the field of $p$-adic harmonic analysis \cite{HS,HS1,HSG1,K,SG,SH,S} and mathematical physics \cite{VSV,VVZ}. Moreover, $p$-adic analysis has tremendous applications in the likes of spring glass theory \cite{ABK}, string theory \cite{VV}, quantum mechanics \cite{VVZ} and quantum gravity \cite{ADFV,BF}.

For a fixed prime $p,$ it is always possible to write a nonzero rational number $x$ in the form $x=p^\gamma m/n,$ where $p$ is not divisible by $m,n\in\mathbb{Z},$ and $\gamma$ is an integer. The $p$-adic norm is defined as follows: $|x|_p=\{0\}\cup\{p^{-\gamma}:\gamma\in\mathbb{Z}\}.$ The $p$-adic norm $|\cdot|_p$ fulfills all the properties of a real norm along with a stronger inequality:
 \begin{equation}\label{NE}|x+y|_p\le\max\{|x|_p,|y|_p\}.\end{equation}
The completion of the field of rational number with respect to $|.|_{p}$ leads to the field of $p$-adic numbers $\mathbb{Q}_p$. In \cite{VVZ}, it was shown that any $x\in \mathbb{Q}_p\backslash\{0\}$ can be represented in the series form as: \begin{equation}\label{E3} x=p^\gamma\sum_{j=0}^{\infty}\beta_jp^j,\end{equation}
 where $\beta_j,\gamma\in\mathbb Z, \beta_j\in\frac{\mathbb{Z}}{p\mathbb{Z}_{p}}, \beta_{0}\neq0.$ The convergence of series (\ref{E3}) is followed from $|p^\gamma\beta_{k}p^{j}|_{p}= p^{-\gamma-j}.$

 The space $\mathbb Q_p^n=\mathbb Q_p\times...\times\mathbb Q_p$ consists of points $\mathbf{x}=(x_{1},x_{2},...,x_{n}),$ where $x_i\in\mathbb Q_p,i=1,2,...,n.$  The norm on $\mathbb Q_p^n$ is given by
  \begin{equation*}|\mathbf{x}|_{p}=\max_{1\leq i \leq n}|x_{i}|_{p}.\end{equation*}

The ball $B_{\gamma}(\mathbf{a})$ and the corresponding sphere $S_{\gamma}(\mathbf{a})$ with center at $\mathbf{a} \in \mathbb{Q}_p^n$ and radius $p^{\gamma}$ in non-Archimedean geometry are defined as follows.
$$B_{\gamma}(\mathbf{a})=\{\mathbf{x} \in \mathbb{Q}_p^n:|\mathbf{x}-\mathbf{a}|_{p}\leq p^{\gamma}\}, \
S_{\gamma}(\mathbf{a})=\{\mathbf{x} \in \mathbb{Q}_p^n:|\mathbf{x}-\mathbf{a}|_{p}=p^{\gamma}\}.$$
 When $\mathbf{a}=\mathbf{0},$ we write $B_\gamma(\mathbf{0})=B_\gamma, \ S_\gamma(\mathbf{0})=S_\gamma.$

Since the space $\mathbb{Q}_p^n$ is locally compact commutative group under addition, it cements the fact from the standard analysis that there exists a translation invariant Haar measure $d\mathbf{x}.$ Also, the measure is normalized by
$$\int_{B_{0}}d\mathbf{x}=|B_{0}|_{H}=1,$$
where $|B|_{H}$ represents the Haar measure of a measurable subset $B$ of $\mathbb{Q}_p^n.$ Furthermore, one can easily show  $|B_{\gamma}(\mathbf{a})|=p^{n\gamma}$, $|S_{\gamma}(\mathbf{a})|=p^{n\gamma}(1-p^{-n})$, for any $\mathbf{a}\in \mathbb{Q}_p^n$.\\
On the other hand, weighted Hardy operator in $p$-adic field was defined in \cite{RLL} and is as follows.
$$H^{p}_{\psi}f(x)=\int_{\mathbb{Z}_{p}^*}f(tx)\psi(t)dt,$$
where $\psi$ is a nonnegative function defined on $\mathbb{Z}_{p}^*.$ In the same paper it was shown that $H^{p}_{\psi}$ is bounded on $L^{q}(\mathbb{Q}_{p}^{n}),$ $1\leq q<\infty$ and $BMO(\mathbb{Q}_{p}^{n}).$ The corresponding norms were also obtained. In 2013, authors in \cite{CD} acquired the boundedness of $H^{p}_{\psi}$ on $p$-adic Herz type spaces. They also established the operator norm. Furthermore, weighted $p$-adic Hardy operator and its commutator on $p$-adic central Morrey spaces were discussed in \cite{WF}.
Clearly, if $\psi\equiv1$ and $n=1,$ we get the $p$-adic Hardy operator on $\mathbb{Q}_{p}$ which is defined by
$$H^{p}f(x)=\frac{1}{|x|_{p}}\int_{|t|_{p}\leq|x|_{p}}f(t)dt,\quad x\neq0.$$
Whereas, for a nonnegative measurable function $f$ on $\mathbb{Q}_p^n$, Fu et al. \cite{FWL} defined the $n$-dimensional $p$-adic Hardy operator as follows.
$$H^{p}f(\mathbf{x})=\frac{1}{|\mathbf{x}|_{p}^{n}}\int_{|\mathbf{t}|_{p}\leq|\mathbf{x}|_{p}}f(\mathbf{t})d\mathbf{t},\quad \mathbf{x}\in\mathbb{Q}_{p}^{n}\setminus\{0\}.$$
In recent times, boundedness of $p$-adic Hardy type operators has attracted many researchers see, for instance \cite{FWL,GZ,LZ,W,WMF}.

Next, we shift our attention towards the multilinear case of an operator. It is worhwhile mentioning here that multilinear operators are widely studied in the past. Multilinear operators are not only generalization of linear one but also occur naturally time and again in analysis. For better understanding of multilinear operators we refer articles including \cite{FLL,FGLY,GL}.
The weighted multilinear Hardy operator was defined by Fu et al.\cite{FGLY} and is follows.
$$H_{\psi}^{m}(f_{1},\cdot\cdot\cdot,f_{m})(\mathbf{x})=\int_{0}^{1}\cdot\cdot\cdot\int_{0}^{1}f_{1}(\mathbf{t}_{1}\mathbf{x})\cdot\cdot\cdot f_{m}(\mathbf{t}_{m}\mathbf{x})\psi(\mathbf{t}_{1},\cdot\cdot\cdot,\mathbf{t}_{m})d\mathbf{t}_{1}\cdot\cdot\cdot d\mathbf{t}_{m},\quad \mathbf{x}\in\mathbb{R}^{n},$$ where $\psi$ is a nonnegative function defined on $[0,1]\times\cdot\cdot\cdot\times[0,1].$ In the same paper, authors studied the boundedness of the very operator on the product of Lebesgue spaces and central Morrey spaces. Moreover in \cite{GFM}, it was also shown that weighted multilinear Hardy operators are bounded on the product of Herz type spaces.
On the other hand, in the $p$-adic harmonic analysis multilinear operators has played a vital role, for example \cite{LZ,LZ1} and the reference therein. The weighted multilinear Hardy operator in the $p$-adic fields was defined in \cite{LZ1} and is as follows. Let $\mathbf{x}\in\mathbb{Q}_{p}^{n},$ $m\in\mathbb{N}$ and $f_{1},\cdot\cdot\cdot,f_{m}$ be nonnegative measurable functions on $\mathbb{Q}_{p}^{n},$ then
$$H^{p,m}_{\psi}(f_{1},\cdot\cdot\cdot,f_{m})(\mathbf{x})=\int_{\mathbb{Z}_{p}^*}\cdot\cdot\cdot\int_{\mathbb{Z}_{p}^*}f_{1}(\mathbf{t}_{1}\mathbf{x})\cdot\cdot\cdot f_{m}(\mathbf{t}_{m}\mathbf{x})\psi(\mathbf{t}_{1},\cdot\cdot\cdot,\mathbf{t}_{m})d\mathbf{t}_{1}\cdot\cdot\cdot d\mathbf{t}_{m},$$
where $\psi$ is a nonnegative measurable function on $\mathbb{Z}_{p}^*\times\cdot\cdot\cdot\times\mathbb{Z}_{p}^*.$ Authors obtained the sharp bounds of weighted multilinear $p$-adic Hardy operator on the product of Lebesgue spaces and Morrey type spaces.

We denote $H^{p,m,*}_{\psi},$ the dual operator of weighted multilinear $p$-adic Hardy operator and is defined as follows. Let $\mathbf{x}\in\mathbb{Q}_{p}^{n},$ $m\in\mathbb{N}$ and $f_{1},\cdot\cdot\cdot,f_{m}$ be nonnegative measurable functions on $\mathbb{Q}_{p}^{n},$ then we have:
$$H^{p,m}_{\psi}(f_{1},\cdot\cdot\cdot,f_{m})(\mathbf{x})=\int_{\mathbb{Z}_{p}^*}\cdot\cdot\cdot\int_{\mathbb{Z}_{p}^*}f_{1}\bigg(\frac{\mathbf{x}}{|\mathbf{t}_{1}|_{p}}\bigg)\cdot\cdot\cdot f_{m}\bigg(\frac{x}{|\mathbf{t}_{m}|_{p}}\bigg)(|\mathbf{t}_{1}|_{p})^{-n}\cdot\cdot\cdot(|\mathbf{t}_{m}|_{p})^{-n}$$$$\quad\times\psi(\mathbf{t}_{1},\cdot\cdot\cdot,\mathbf{t}_{m})d\mathbf{t}_{1}\cdot\cdot\cdot d\mathbf{t}_{m},$$
where $\psi$ is a nonnegative measurable function on $\mathbb{Z}_{p}^*\times\cdot\cdot\cdot\times\mathbb{Z}_{p}^*.$

Motivated from above results, the purpose of this article is to establish a sufficient and necessary condition of weighted multilinear $p$-adic Hardy operators on the product of Herz spaces and Morrey-Herz spaces. Before moving to our main results, let us specify that $\chi_{k}$ is the characteristic function of a sphere $S_{k}.$ Also we recall the definition of homogeneous $p$-adic Herz spaces and homogeneous $p$-adic Morrey-Herz spaces.
\begin{definition}\cite{FWL} Suppose $\alpha\in\mathbb{R},$ $0<q<\infty$ and $0<r<\infty.$ The homogeneous $p$-adic Herz space $\dot{K}^{\alpha,r}_{q}(\mathbb{Q}_{p}^{n})$ is defined by
$$K^{\alpha,r}_{q}(\mathbb{Q}_{p}^{n})=\{f\in L^{q}(\mathbb{Q}_{p}^{n}):\|f\|_{K^{\alpha,r}_{q}(\mathbb{Q}_{p}^{n})}<\infty\},$$
where $$\|f\|_{K^{\alpha,r}_{q}(\mathbb{Q}_{p}^{n})}=\bigg(\sum_{k=-\infty}^{\infty}p^{k\alpha r}\|f\chi_{k}\|^{r}_{L^{q}(\mathbb{Q}_{p}^{n})}\bigg)^{1/q}.$$
\end{definition}
Obviously, $\dot{K}^{0,q}_{q}(\mathbb{Q}_{p}^{n})=L^{q}(\mathbb{Q}_{p}^{n})$ and $\dot{K}^{\alpha/q,q}_{q}(\mathbb{Q}_{p}^{n})=L^{q}(|\mathbf{x}|_{p}^{\alpha}).$ \begin{definition}\cite{CD}\label{vv} Let $\alpha\in\mathbb{R},$ $0<r<\infty,$ $0<q<\infty$ and $\lambda\geq0.$ The homogeneous $p$-adic Morrey-Herz space is defined as: $$M\dot{K}^{\alpha,\lambda}_{r,q}(\mathbb{Q}_{p}^{n})=\{f\in L^{q}_{loc}(\mathbb{Q}_{p}^{n}\setminus\{0\}):\|f\|_{M\dot{K}^{\alpha,\lambda}_{r,q}(\mathbb{Q}_{p}^{n})}<\infty\},$$
where $$\|f\|_{M\dot{K}^{\alpha,\lambda}_{r,q}(\mathbb{Q}_{p}^{n})}=\sup_{k_{0}\in\mathbb{Z}}p^{-k_{0}\lambda}\bigg(\sum_{k=-\infty}^{\infty}p^{k\alpha r}\|f\chi_{k}\|^{r}_{L^{q}(\mathbb{Q}_{p}^{n})}\bigg)^{1/r}.$$
\end{definition}
It is eminent that $M\dot{K}^{\alpha,0}_{r,q}(\mathbb{Q}_{p}^{n})=\dot{K}^{\alpha,r}_{q}(\mathbb{Q}_{p}^{n})$ and $M\dot{K}^{\alpha/q,0}_{q,q}(\mathbb{Q}_{p}^{n})=L^{q}(|\mathbf{x}|_{p}^{\alpha}).$
\section{Boundedness of $H^{p,m}_{\psi}$ and $H^{p,m,*}_{\psi}$ on the product of $p$-adic Herz spaces}
In the current section we guarantee the boundedness of weighted multilinear Hardy operator and its adjoint on the product of Herz spaces. We also obtain the corresponding norms of both operators. The main results of the section are as follows.
\begin{theorem}\label{T111}Let $\alpha,\alpha_{1},\alpha_{2},...,\alpha_{m}$ be any arbitrary real numbers, $1<p,p_{1},....,p_{m},q,q_{1},...,q_{m}<\infty$ and let also $\alpha_{1}+\alpha_{2}+...+\alpha_{m}=\alpha,$ $\frac{1}{p_{1}}+\frac{1}{p_{2}}+...+\frac{1}{p_{m}}=\frac{1}{p},$ $\frac{1}{q_{1}}+\frac{1}{q_{2}}+...+\frac{1}{q_{m}}=\frac{1}{q},$ then $H^{p,m}_{\psi}$ is bounded from $\dot{K}_{q_{1}}^{\alpha_{1},p_{1}}(\mathbb{Q}_p^n)\times...\times\dot{K}_{q_{m}}^{\alpha_{m},p_{m}}(\mathbb{Q}_p^n)$ to $\dot{K}_{q}^{\alpha,p}(\mathbb{Q}_p^n)$ if
\begin{eqnarray}\begin{aligned}\label{zxx}\int_{\mathbb{{Z}}_{p}^*}...\int_{\mathbb{{Z}}_{p}^*}|\mathbf{t}_{1}|_{p}^{-(\alpha_{1}+n/q_{1})}...|\mathbf{t}_{m}|_{p}^{-(\alpha_{m}+n/q_{m})}\psi(\mathbf{t}_{1},...,\mathbf{t}_{m})d\mathbf{t}_{1}...d\mathbf{t}_{m}<\infty.\end{aligned}\end{eqnarray}
Conversely,if $q_{1},q_{2},...q_{m}=mq,$ $p_{1},p_{2},...,p_{m}=mp$ and $H^{p,m}_{\psi}$ is bounded from $\dot{K}_{q_{1}}^{\alpha_{1},p_{1}}(\mathbb{Q}_p^n)\times...\times\dot{K}_{q_{m}}^{\alpha_{m},p_{m}}(\mathbb{Q}_p^n)$ to $\dot{K}_{q}^{\alpha,p}(\mathbb{Q}_p^n)$ then (\ref{zxx}) holds. Furthermore,
\begin{eqnarray}\begin{aligned}[b]\|H^{p,m}_{\psi}\|_{\dot{K}_{q_{1}}^{\alpha_{1},p_{1}}(\mathbb{Q}_p^n)\times...\times\dot{K}_{q_{m}}^{\alpha_{m},p_{m}}(\mathbb{Q}_p^n)\rightarrow\dot{K}_{q}^{\alpha,p}(\mathbb{Q}_p^n)}=&\int_{\mathbb{{Z}}_{p}^*}...\int_{\mathbb{{Z}}_{p}^*}|\mathbf{t}_{1}|_{p}^{-(\alpha_{1}+n/q_{1})}...|\mathbf{t}_{m}|_{p}^{-(\alpha_{m}+n/q_{m})}\\&\quad\times\psi(\mathbf{t}_{1},...,\mathbf{t}_{m})d\mathbf{t}_{1}...d\mathbf{t}_{m}.\end{aligned}\end{eqnarray}
\end{theorem}
In a similar manner, we have a result for an operator $H^{p,m,*}_{\psi}.$
\begin{theorem}\label{T1111}Let $\alpha,\alpha_{1},\alpha_{2},...,\alpha_{m}$ be any arbitrary real numbers, $1<p,p_{1},....,p_{m},q,q_{1},...,q_{m}<\infty$ and let also $\alpha_{1}+\alpha_{2}+...+\alpha_{m}=\alpha,$ $\frac{1}{p_{1}}+\frac{1}{p_{2}}+...+\frac{1}{p_{m}}=\frac{1}{p},$ $\frac{1}{q_{1}}+\frac{1}{q_{2}}+...+\frac{1}{q_{m}}=\frac{1}{q},$ then $H^{p,m,*}_{\psi}$ is bounded from $\dot{K}_{q_{1}}^{\alpha_{1},p_{1}}(\mathbb{Q}_p^n)\times...\times\dot{K}_{q_{m}}^{\alpha_{m},p_{m}}(\mathbb{Q}_p^n)$ to $\dot{K}_{q}^{\alpha,p}(\mathbb{Q}_p^n)$ if
\begin{eqnarray}\begin{aligned}\label{zx}\int_{\mathbb{{Z}}_{p}^*}...\int_{\mathbb{{Z}}_{p}^*}|\mathbf{t}_{1}|_{p}^{\alpha_{1}-n(1-1/q_{1})}...|\mathbf{t}_{m}|_{p}^{\alpha_{m}-n(1-1/q_{m})}\psi(\mathbf{t}_{1},...,\mathbf{t}_{m})d\mathbf{t}_{1}...d\mathbf{t}_{m}<\infty.\end{aligned}\end{eqnarray}
Conversely,if $q_{1},q_{2},...q_{m}=mq,$ $p_{1},p_{2},...,p_{m}=mp$ and $H^{p,m}_{\psi}$ is bounded from $\dot{K}_{q_{1}}^{\alpha_{1},p_{1}}(\mathbb{Q}_p^n)\times...\times\dot{K}_{q_{m}}^{\alpha_{m},p_{m}}(\mathbb{Q}_p^n)$ to $\dot{K}_{q}^{\alpha,p}(\mathbb{Q}_p^n)$ then (\ref{zx}) holds. Furthermore,
\begin{eqnarray}\begin{aligned}[b]\|H^{p,m,*}_{\psi}\|_{\dot{K}_{q_{1}}^{\alpha_{1},p_{1}}(\mathbb{Q}_p^n)\times...\times\dot{K}_{q_{m}}^{\alpha_{m},p_{m}}(\mathbb{Q}_p^n)\rightarrow\dot{K}_{q}^{\alpha,p}(\mathbb{Q}_p^n)}=&\int_{\mathbb{{Z}}_{p}^*}...\int_{\mathbb{{Z}}_{p}^*}|\mathbf{t}_{1}|_{p}^{\alpha_{1}-n(1-1/q_{1})}...|\mathbf{t}_{m}|_{p}^{\alpha_{m}-n(1-1/q_{m})}\\&\quad\times\psi(\mathbf{t}_{1},...,\mathbf{t}_{m})d\mathbf{t}_{1}...d\mathbf{t}_{m}.\end{aligned}\end{eqnarray}
\end{theorem}
The proof of Theorem (\ref{T1111}) can be obtained similarly to  the proof of Theorem (\ref{T111}). So, we only prove Theorem (\ref{T111})

\noindent\textit{ \textbf{Proof of Theorem \ref{T111}:}}
We will prove the Theorem only for $m=2$ which will work for all $m\in\mathbb{N}.$ Since $1/q=1/q_{1}+1/q_{2},$ applying Minkowski's inequality and H\"{o}lder's inequality to have:
\begin{eqnarray*}\begin{aligned}\|H^{p,2}_{\psi}(f_{1},f_{2})\chi_{k}\|_{L^{q}(\mathbb{Q}_p^n)}=&\bigg(\int_{S_{k}}\bigg|\int\int_{\mathbb{Z}_{p}^*}f_{1}(\mathbf{t}_{1}\mathbf{x})f_{2}(\mathbf{t}_{2}\mathbf{x})\psi(\mathbf{t}_{1},\mathbf{t}_{2})d\mathbf{t}_{1}d\mathbf{t}_{2}\bigg|^{q}\bigg)^{1/q}\\
\leq&\int\int_{\mathbb{Z}_{p}^*}\bigg(\int_{S_{k}}\bigg|f_{1}(\mathbf{t}_{1}\mathbf{x})f_{2}(\mathbf{t}_{2}\mathbf{x})\bigg|^{q}\bigg)^{1/q}\psi(\mathbf{t}_{1},\mathbf{t}_{2})d\mathbf{t}_{1}d\mathbf{t}_{2}\\
\leq&\int\int_{\mathbb{Z}_{p}^*}\bigg(\int_{S_{k}}|f_{1}(\mathbf{t}_{1}\mathbf{x})|^{q_{1}}d\mathbf{x}\bigg)^{1/q_{1}}\bigg(\int_{S_{k}}|f_{2}(\mathbf{t}_{2}\mathbf{x})|^{q_{2}}d\mathbf{x}\bigg)^{1/q_{2}}\\\quad&\times\psi(\mathbf{t}_{1},\mathbf{t}_{2})d\mathbf{t}_{1}d\mathbf{t}_{2}\\
=&\int\int_{\mathbb{Z}_{p}^*}\bigg(\int_{\mathbf{t}_{1}S_{k}}|f_{1}(\mathbf{x})|^{q_{1}}d\mathbf{x}\bigg)^{1/q_{1}}\bigg(\int_{\mathbf{t}_{2}S_{k}}|f_{2}(\mathbf{x})|^{q_{2}}d\mathbf{x}\bigg)^{1/q_{2}}\\\quad&\times|\mathbf{t}_{1}|_{p}^{-n/q_{1}}|\mathbf{t}_{2}|_{p}^{-n/q_{2}}\psi(\mathbf{t}_{1},\mathbf{t}_{2})d\mathbf{t}_{1}d\mathbf{t}_{2}.
\end{aligned}\end{eqnarray*}
Now for each $\mathbf{t}_{1},\mathbf{t}_{2}\in\mathbb{Z}_{p}^{*},$ there exists non-negative integers $m,l$ such that $|\mathbf{t}_{1}|_{p}=p^{-m}$ and $|\mathbf{t}_{2}|_{p}=p^{-l}.$ Therefore, we easily have:
\begin{eqnarray*}\begin{aligned}\|H^{p,2}_{\psi}(f_{1},f_{2})\chi_{k}\|_{L^{q}(\mathbb{Q}_p^n)}\leq&\int\int_{\mathbb{Z}_{p}^*}\bigg(\|f_{1}\chi_{k-m}\|_{L^{q_{1}}(\mathbb{Q}_p^n)}.\|f_{2}\chi_{k-l}\|_{L^{q_{2}}(\mathbb{Q}_p^n)}\bigg)\\\quad&\times|\mathbf{t}_{1}|_{p}^{-n/q_{1}}|\mathbf{t}_{2}|_{p}^{-n/q_{2}}\psi(\mathbf{t}_{1},\mathbf{t}_{2})d\mathbf{t}_{1}d\mathbf{t}_{2}.
\end{aligned}\end{eqnarray*}
Hence, by means of Minkowski's inequality and H\"{o}lder's inequality together with $1/p=1/p_{1}+1/p_{2}$ and $\alpha=\alpha_{1}+\alpha_{2},$  we get:
\begin{eqnarray*}\begin{aligned}[b]\label{cc}\|H^{p,2}_{\psi}(f_{1},f_{2})\|_{\dot{K}^{\alpha,p}_{q}(\mathbb{Q}_p^n)}=&\bigg(\sum_{k=-\infty}^{\infty}p^{k\alpha p}\|H^{p,2}_{\psi}(f_{1},f_{2})\chi_{k}\|^{p}_{L^{q}(\mathbb{Q}_p^n)}\bigg)^{1/p}\\
\end{aligned}\end{eqnarray*}
\begin{eqnarray}\begin{aligned}[b]\label{cc}
\leq&\bigg(\sum_{k=-\infty}^{\infty}p^{k\alpha p}\bigg(\int\int_{\mathbb{Z}_{p}^*}\bigg(\|f_{1}\chi_{k-m}\|_{L^{q_{1}}(\mathbb{Q}_p^n)}.\|f_{2}\chi_{k-l}\|_{L^{q_{2}}(\mathbb{Q}_p^n)}\bigg)\\\quad&\times|\mathbf{t}_{1}|_{p}^{-n/q_{1}}|\mathbf{t}_{2}|_{p}^{-n/q_{2}}\psi(\mathbf{t}_{1},\mathbf{t}_{2})d\mathbf{t}_{1}d\mathbf{t}_{2}\bigg)^{p}\bigg)^{1/p}\\
\leq&\int\int_{\mathbb{Z}_{p}^*}\bigg(\sum_{k=-\infty}^{\infty}p^{k\alpha p}\bigg(\|f_{1}\chi_{k-m}\|_{L^{q_{1}}(\mathbb{Q}_p^n)}.\|f_{2}\chi_{k-l}\|_{L^{q_{2}}(\mathbb{Q}_p^n)}\bigg)^{p}\bigg)^{1/p}\\\quad&\times|\mathbf{t}_{1}|_{p}^{-n/q_{1}}|\mathbf{t}_{2}|_{p}^{-n/q_{2}}\psi(\mathbf{t}_{1},\mathbf{t}_{2})d\mathbf{t}_{1}d\mathbf{t}_{2}\\
\leq&\int\int_{\mathbb{Z}_{p}^*}\bigg(\sum_{k=-\infty}^{\infty}p^{k\alpha_{1} p_{1}}\|f_{1}\chi_{k-m}\|^{p_{1}}_{L^{q_{1}}(\mathbb{Q}_p^n)}\bigg)^{1/p_{1}}\\\quad&\times\bigg(\sum_{k=-\infty}^{\infty}p^{k\alpha_{2} p_{2}}\|f_{2}\chi_{k-l}\|^{p_{2}}_{L^{q_{2}}(\mathbb{Q}_p^n)}\bigg)^{1/p_{2}}|\mathbf{t}_{1}|_{p}^{-n/q_{1}}|\mathbf{t}_{2}|_{p}^{-n/q_{2}}\psi(\mathbf{t}_{1},\mathbf{t}_{2})d\mathbf{t}_{1}d\mathbf{t}_{2}\\
\leq&\int\int_{\mathbb{Z}_{p}^*}\bigg(\sum_{k=-\infty}^{\infty}p^{k\alpha_{1} p_{1}}\|f_{1}\chi_{k}\|^{p_{1}}_{L^{q_{1}}(\mathbb{Q}_p^n)}\bigg)^{1/p_{1}}\\\quad&\times\bigg(\sum_{k=-\infty}^{\infty}p^{k\alpha_{2} p_{2}}\|f_{2}\chi_{k}\|^{p_{2}}_{L^{q_{2}}(\mathbb{Q}_p^n)}\bigg)^{1/p_{2}}|\mathbf{t}_{1}|_{p}^{-(\alpha_{1}+n/q_{1})}|\mathbf{t}_{2}|_{p}^{-(\alpha_{2}+n/q_{2})}\\\quad&\times\psi(\mathbf{t}_{1},\mathbf{t}_{2})d\mathbf{t}_{1}d\mathbf{t}_{2}\\
\leq&\|f_{1}\|_{\dot{K}_{q_{1}}^{\alpha_{1},p_{1}}(\mathbb{Q}_p^n)}\|f_{2}\|_{\dot{K}_{q_{2}}^{\alpha_{2},p_{2}}(\mathbb{Q}_p^n)}\int\int_{\mathbb{Z}_{p}^*}|\mathbf{t}_{1}|_{p}^{-(\alpha_{1}+n/q_{1})}|\mathbf{t}_{2}|_{p}^{-(\alpha_{2}+n/q_{2})}\\\quad&\times\psi(\mathbf{t}_{1},\mathbf{t}_{2})d\mathbf{t}_{1}d\mathbf{t}_{2}.
\end{aligned}\end{eqnarray}
Hence,
\begin{eqnarray}\begin{aligned}[b]\label{mi}\|H^{p,2}_{\psi}(f_{1},f_{2})\chi_{k}\|_{\dot{K}^{\alpha_{1},p_{1}}_{q_{1}}(\mathbb{Q}_p^n)\times\dot{K}^{\alpha_{2},p_{2}}_{q_{2}}(\mathbb{Q}_p^n)\rightarrow\dot{K}^{\alpha,p}_{q}(\mathbb{Q}_p^n)}\leq&\int\int_{\mathbb{Z}_{p}^*}|\mathbf{t}_{1}|_{p}^{-(\alpha_{1}+n/q_{1})}|\mathbf{t}_{2}|_{p}^{-(\alpha_{2}+n/q_{2})}\\\quad&\times\psi(\mathbf{t}_{1},\mathbf{t}_{2})d\mathbf{t}_{1}d\mathbf{t}_{2}.
\end{aligned}\end{eqnarray}
From (\ref{cc}), first part of Theorem is followed.\\
Conversely, let $H^{p,2}_{\psi}$ be bounded from $\dot{K}_{q_{1}}^{\alpha_{1},p_{1}}(\mathbb{Q}_p^n)\times\dot{K}_{q_{2}}^{\alpha_{2},p_{2}}(\mathbb{Q}_p^n)$ to $\dot{K}_{q}^{\alpha,p}(\mathbb{Q}_p^n).$ For $0<\epsilon<1,$ we let
\begin{eqnarray*}\begin{aligned}
 f_{1}(\mathbf{x})&={\begin{cases} 0  & \text{if } |\mathbf{x}|_{p}<1,\\
|\mathbf{x}|_{p}^{-\alpha_{1}-(n/q_{1})-\epsilon}  & \text{if } |\mathbf{x}|_{p}\geq1,\end{cases}}
\end{aligned}\end{eqnarray*}
\begin{eqnarray*}\begin{aligned}
 f_{2}(\mathbf{x})&={\begin{cases} 0  & \text{if } |\mathbf{x}|_{p}<1,\\
|\mathbf{x}|_{p}^{-\alpha_{2}-(n/q_{2})-\epsilon}  & \text{if } |\mathbf{x}|_{p}\geq1.\end{cases}}
\end{aligned}\end{eqnarray*}
It is quite evident that $f_{1}\chi_{k}=f_{2}\chi_{k}=0$ for $k<0.$ Our interest lies only for $k\geq0.$ So, we proceed as follows:
\begin{eqnarray*}\begin{aligned}\|f_{1}\chi_{k}\|_{L^{q_{1}}(\mathbb{Q}_p^n)}=\bigg(\int_{S_{k}}|\mathbf{x}|_{p}^{-(\alpha_{1}+(n/q_{1})+\epsilon)q_{1}}\bigg)^{1/q_{1}}=(1-p^{-n})^{1/q_{1}}p^{-k(\alpha_{1}+\epsilon)}.
\end{aligned}\end{eqnarray*}
In a similar fashion, we have \begin{eqnarray*}\begin{aligned}\|f_{2}\chi_{k}\|_{L^{q_{2}}(\mathbb{Q}_p^n)}=(1-p^{-n})^{1/q_{2}}p^{-k(\alpha_{2}+\epsilon)}.
\end{aligned}\end{eqnarray*}
Hence,
\begin{eqnarray*}\begin{aligned}\|f_{1}\|_{\dot{K}_{q_{1}}^{\alpha_{1},p_{1}}(\mathbb{Q}_p^n)}=&\bigg(\sum_{k=-\infty}^{\infty}p^{k\alpha_{1}p_{1}}\|f_{1}\chi_{k}\|^{p_{1}}_{L^{q_{1}}(\mathbb{Q}_p^n)}\bigg)^{1/p_{1}}\\
=&(1-p^{-n})^{1/q_{1}}\bigg(\sum_{k=0}^{\infty}p^{k\alpha_{1}p_{1}}p^{-k(\alpha_{1}+\epsilon)p_{1}}\bigg)^{1/p_{1}}\\
=&(1-p^{-n})^{1/q_{1}}\bigg(\sum_{k=0}^{\infty}p^{-kp_{1}\epsilon}\bigg)^{1/p_{1}}\\
=&(1-p^{-n})^{1/q_{1}}\frac{p^{\epsilon}}{p^{p_{1}\epsilon}-1}.
\end{aligned}\end{eqnarray*}
Similarly,
\begin{eqnarray*}\begin{aligned}\|f_{2}\|_{\dot{K}_{q_{2}}^{\alpha_{2},p_{2}}(\mathbb{Q}_p^n)}=&(1-p^{-n})^{1/q_{2}}\frac{p^{\epsilon}}{p^{p_{2}\epsilon}-1}.
\end{aligned}\end{eqnarray*}
It is obvious to see that when $|\mathbf{x}|_{p}<1$ then $H^{p,2}_{\psi}(f_{1},f_{2})=0.$ So, we are down to the case $|\mathbf{x}|_{p}\geq1.$ We will evaluate the case as:
\begin{eqnarray*}\begin{aligned}H^{p,2}_{\psi}(f_{1},f_{2})=&\int\int_{|\mathbf{x}|_{p}^{-1}\leq|\mathbf{t}|_{p}\leq1}f_{1}(\mathbf{t}_{1}\mathbf{x})(\mathbf{t}_{2}\mathbf{x})\psi(\mathbf{t}_{1},\mathbf{t}_{2})d\mathbf{t}_{1}d\mathbf{t}_{2}\\
=&|\mathbf{x}|_{p}^{-(\alpha_{1}+n/q+2\epsilon)}\int\int_{|\mathbf{x}|_{p}^{-1}\leq|\mathbf{t}|_{p}\leq1}|\mathbf{t}_{1}|_{p}^{-\alpha_{1}-n/q_{1}-\epsilon}|\mathbf{t}_{2}|_{p}^{-\alpha_{1}-n/q_{2}-\epsilon}\psi(\mathbf{t}_{1},\mathbf{t}_{2})d\mathbf{t}_{1}d\mathbf{t}_{2}.
\end{aligned}\end{eqnarray*}
Now for $k\leq0,$ we have $H^{p,2}_{\psi}(f_{1},f_{2})\chi_{k}=0.$ So, the interest is left only for $k\geq0.$ By the definition of Herz space, we have:
\begin{eqnarray*}\begin{aligned}\|H^{p,2}_{\psi}(f_{1},f_{2})\chi_{k}\|^{p}_{\dot{K}^{\alpha,p}_{q}(\mathbb{Q}_p^n)}=&\sum_{k=0}^{\infty}p^{k\alpha p}\|H^{p,2}_{\psi}(f_{1},f_{2})\chi_{k}\|^{p}_{L^{q}(\mathbb{Q}_p^n)}\\
=&\sum_{k=0}^{\infty}p^{k\alpha p}\bigg\{\int_{S_{k}}(|\mathbf{x}|_{p}^{-(\alpha+n/q+2\epsilon)}\\&\quad\times\int\int_{p^{-k}\leq|\mathbf{t}|_{p}\leq1}|\mathbf{t}_{1}|_{p}^{-\alpha_{1}-n/q_{1}-\epsilon}|\mathbf{t}_{2}|_{p}^{-\alpha_{2}-n/q_{2}-\epsilon}\psi(\mathbf{t}_{1},\mathbf{t}_{2})d\mathbf{t}_{1}d\mathbf{t}_{2})^{q}d\mathbf{x}\bigg\}^{p/q}\\
\end{aligned}\end{eqnarray*}
\begin{eqnarray*}\begin{aligned}
=&(1-p^{-n})^{p/q}\sum_{k=0}^{\infty}p^{k\alpha p}(p^{-k(\alpha+2\epsilon)p})\\&\quad\times\bigg(\int\int_{p^{-k}\leq|\mathbf{t}|_{p}\leq1}|\mathbf{t}_{1}|_{p}^{-\alpha_{1}-n/q_{1}-\epsilon}|\mathbf{t}_{2}|_{p}^{-\alpha_{2}-n/q_{2}-\epsilon}\psi(\mathbf{t}_{1},\mathbf{t}_{2})d\mathbf{t}_{1}d\mathbf{t}_{2}\bigg)^{p}.
\end{aligned}\end{eqnarray*}
Now for any $l\leq k$, we get
\begin{eqnarray*}\begin{aligned}\|H^{p,2}_{\psi}(f_{1},f_{2})\chi_{k}\|_{\dot{K}^{\alpha,p}_{q}(\mathbb{Q}_p^n)}\geq&(1-p^{-n})^{1/q}\bigg(\sum_{k=l}^{\infty}p^{-2k\epsilon p}\bigg)^{1/p}\\&\quad\times\bigg(\int\int_{p^{-l}\leq|\mathbf{t}|_{p}\leq1}|\mathbf{t}_{1}|_{p}^{-\alpha_{1}-n/q_{1}-\epsilon}|\mathbf{t}_{2}|_{p}^{-\alpha_{2}-n/q_{2}-\epsilon}\psi(\mathbf{t}_{1},\mathbf{t}_{2})d\mathbf{t}_{1}d\mathbf{t}_{2}\bigg)\\
=&(1-p^{-n})^{1/q}\bigg(\sum_{k=0}^{\infty}p^{-2k\epsilon p}\bigg)^{1/p}\\&\quad\times\bigg(p^{-2l\epsilon}\int\int_{p^{-l}\leq|\mathbf{t}|_{p}\leq1}|\mathbf{t}_{1}|_{p}^{-\alpha_{1}-n/q_{1}-\epsilon}|\mathbf{t}_{2}|_{p}^{-\alpha_{2}-n/q_{2}-\epsilon}\psi(\mathbf{t}_{1},\mathbf{t}_{2})d\mathbf{t}_{1}d\mathbf{t}_{2}\bigg)\\
=&(1-p^{-n})^{1/q}\frac{p^{2\epsilon}}{(p^{2\epsilon p}-1)^{1/p}}\\&\quad\times\bigg(p^{-2l\epsilon}\int\int_{p^{-l}\leq|\mathbf{t}|_{p}\leq1}|\mathbf{t}_{1}|_{p}^{-\alpha_{1}-n/q_{1}-\epsilon}|\mathbf{t}_{2}|_{p}^{-\alpha_{2}-n/q_{2}-\epsilon}\psi(\mathbf{t}_{1},\mathbf{t}_{2})d\mathbf{t}_{1}d\mathbf{t}_{2}\bigg).
\end{aligned}\end{eqnarray*}
Since $q_{1}=q_{2}=2q,$ $1/p=1/p_{1}+1/p_{2}$ and $p_{1}=p_{2}=2p,$ we have:
\begin{eqnarray*}\begin{aligned}\|H^{p,2}_{\psi}(f_{1},f_{2})\chi_{k}\|_{\dot{K}^{\alpha,p}_{q}(\mathbb{Q}_p^n)}\geq&\|f_{1}\|_{\dot{K}_{q_{1}}^{\alpha_{1},p_{1}}(\mathbb{Q}_p^n)}\|f_{2}\|_{\dot{K}_{q_{2}}^{\alpha_{2},p_{2}}(\mathbb{Q}_p^n)}\\&\quad\times\bigg(p^{-2l\epsilon}\int\int_{p^{-l}\leq|\mathbf{t}|_{p}\leq1}|\mathbf{t}_{1}|_{p}^{-\alpha_{1}-n/q_{1}-\epsilon}|\mathbf{t}_{2}|_{p}^{-\alpha_{2}-n/q_{2}-\epsilon}\psi(\mathbf{t}_{1},\mathbf{t}_{2})d\mathbf{t}_{1}d\mathbf{t}_{2}\bigg).
\end{aligned}\end{eqnarray*}
We take $\epsilon=p^{-l},$ $l=0,1,2,3,...$, then letting $l\rightarrow\infty,$ we have $\epsilon\rightarrow0.$ Ultimately, we get
\begin{eqnarray}\begin{aligned}[b]\label{ra}\|H^{p,2}_{\psi}(f_{1},f_{2})\chi_{k}\|_{\dot{K}^{\alpha_{1},p_{1}}_{q_{1}}(\mathbb{Q}_p^n)\times\dot{K}^{\alpha_{2},p_{2}}_{q_{2}}(\mathbb{Q}_p^n)\rightarrow\dot{K}^{\alpha,p}_{q}(\mathbb{Q}_p^n)}\geq&\int\int_{\mathbb{Z}_{p}^*}|\mathbf{t}_{1}|_{p}^{-(\alpha_{1}+n/q_{1})}|\mathbf{t}_{2}|_{p}^{-(\alpha_{2}+n/q_{2})}\\\quad&\times\psi(\mathbf{t}_{1},\mathbf{t}_{2})d\mathbf{t}_{1}d\mathbf{t}_{2}.
\end{aligned}\end{eqnarray}
From (\ref{mi}) and (\ref{ra}), we get the required proof.
\section{Boundedness of $H^{p,m}_{\psi}$ and $H^{p,m,*}_{\psi}$ on the product of $p$-adic Morrey-Herz spaces}
The present section addresses the boundedness of $H^{p,m}_{\psi}$ and $H^{p,m,*}_{\psi}$ on the product of $p$-adic Morrey-Herz spaces. The corresponding norms are also ensured. The main results are down under.
\begin{theorem}\label{T101}Let $\alpha,\alpha_{1},\alpha_{2},...,\alpha_{m}$ be any arbitrary real numbers, $1<p,p_{1},....,p_{m},q,q_{1},...,q_{m}<\infty$ and let also $\alpha_{1}+\alpha_{2}+...+\alpha_{m}=\alpha,$ $\lambda,\lambda_{1},...,\lambda_{m}>0,$ $\frac{1}{p_{1}}+\frac{1}{p_{2}}+...+\frac{1}{p_{m}}=\frac{1}{p},$ $\frac{1}{q_{1}}+\frac{1}{q_{2}}+...+\frac{1}{q_{m}}=\frac{1}{q},$ $\lambda_{1}+...+\lambda_{m}=\lambda,$ then $H^{p,m}_{\psi}$ is bounded from $M\dot{K}^{\alpha_{1},\lambda_{1}}_{p_{1},q_{1}}(\mathbb{Q}_p^n)\times...\times M\dot{K}^{\alpha_{m},\lambda_{m}}_{p_{m},q_{m}}(\mathbb{Q}_p^n)$ to $M\dot{K}^{\alpha,\lambda}_{p,q}(\mathbb{Q}_p^n)$ if
\begin{eqnarray}\begin{aligned}\label{nx}\int_{\mathbb{{Z}}_{p}^*}...\int_{\mathbb{{Z}}_{p}^*}|\mathbf{t}_{1}|_{p}^{-(\alpha_{1}+n/q_{1}-\lambda_{1})}...|\mathbf{t}_{m}|_{p}^{-(\alpha_{m}+n/q_{m}-\lambda_{m})}\psi(\mathbf{t}_{1},...,\mathbf{t}_{m})d\mathbf{t}_{1}...d\mathbf{t}_{m}<\infty.\end{aligned}\end{eqnarray}
Conversely,if $q_{1},q_{2},...q_{m}=mq,$ $p_{1},p_{2},...,p_{m}=mp,$ $\alpha_{1}=...=\alpha_{m}=(1/m)\alpha,$ $\lambda_{1}=...=\lambda_{m}=(1/m)\lambda$ and $H^{p,m}_{\psi}$ is bounded from $M\dot{K}^{\alpha_{1},\lambda_{1}}_{p_{1},q_{1}}(\mathbb{Q}_p^n)\times...\times M\dot{K}^{\alpha_{m},\lambda_{m}}_{p_{m},q_{m}}(\mathbb{Q}_p^n)$ to $M\dot{K}^{\alpha,\lambda}_{p,q}(\mathbb{Q}_p^n)$ then (\ref{nx}) holds. Furthermore,
\begin{eqnarray}\begin{aligned}[b]\label{vvv}\|H^{p,m}_{\psi}\|_{M\dot{K}^{\alpha_{1},\lambda_{1}}_{p_{1},q_{1}}(\mathbb{Q}_p^n)\times...\times M\dot{K}^{\alpha_{m},\lambda_{m}}_{p_{m},q_{m}}(\mathbb{Q}_p^n)\rightarrow M\dot{K}^{\alpha,\lambda}_{p,q}(\mathbb{Q}_p^n)}=&\int_{\mathbb{{Z}}_{p}^*}...\int_{\mathbb{{Z}}_{p}^*}|\mathbf{t}_{1}|_{p}^{-(\alpha_{1}+n/q_{1}-\lambda_{1})}...|\mathbf{t}_{m}|_{p}^{-(\alpha_{m}+n/q_{m}-\lambda_{m})}\\&\quad\times\psi(\mathbf{t}_{1},...,\mathbf{t}_{m})d\mathbf{t}_{1}...d\mathbf{t}_{m}.\end{aligned}\end{eqnarray}
\end{theorem}
We also have a corresponding result for an operator $H^{p,m,*}_{\psi}.$
\begin{theorem}\label{T11}Let $\alpha,\alpha_{1},\alpha_{2},...,\alpha_{m}$ be any arbitrary real numbers, $1<p,p_{1},....,p_{m},q,q_{1},...,q_{m}<\infty$ and let also $\alpha_{1}+\alpha_{2}+...+\alpha_{m}=\alpha,$ $\lambda,\lambda_{1},...,\lambda_{m}>0,$ $\frac{1}{p_{1}}+\frac{1}{p_{2}}+...+\frac{1}{p_{m}}=\frac{1}{p},$ $\frac{1}{q_{1}}+\frac{1}{q_{2}}+...+\frac{1}{q_{m}}=\frac{1}{q},$ $\lambda_{1}+...+\lambda_{m}=\lambda,$ then $H^{p,m}_{\psi}$ is bounded from $M\dot{K}^{\alpha_{1},\lambda_{1}}_{p_{1},q_{1}}(\mathbb{Q}_p^n)\times...\times M\dot{K}^{\alpha_{m},\lambda_{m}}_{p_{m},q_{m}}(\mathbb{Q}_p^n)$ to $M\dot{K}^{\alpha,\lambda}_{p,q}(\mathbb{Q}_p^n)$ if
\begin{eqnarray}\begin{aligned}\label{nxx}\int_{\mathbb{{Z}}_{p}^*}...\int_{\mathbb{{Z}}_{p}^*}|\mathbf{t}_{1}|_{p}^{\alpha_{1}-\lambda_{1}-n(1-/q_{1})}...|\mathbf{t}_{m}|_{p}^{\alpha_{m}-\lambda_{m}-n(1-1/q_{m})}\psi(\mathbf{t}_{1},...,\mathbf{t}_{m})d\mathbf{t}_{1}...d\mathbf{t}_{m}<\infty.\end{aligned}\end{eqnarray}
Conversely,if $q_{1},q_{2},...q_{m}=mq,$ $p_{1},p_{2},...,p_{m}=mp,$ $\alpha_{1}=...=\alpha_{m}=(1/m)\alpha,$ $\lambda_{1}=...=\lambda_{m}=(1/m)\lambda$ and $H^{p,m}_{\psi}$ is bounded from $M\dot{K}^{\alpha_{1},\lambda_{1}}_{p_{1},q_{1}}(\mathbb{Q}_p^n)\times...\times M\dot{K}^{\alpha_{m},\lambda_{m}}_{p_{m},q_{m}}(\mathbb{Q}_p^n)$ to $M\dot{K}^{\alpha,\lambda}_{p,q}(\mathbb{Q}_p^n)$ then (\ref{nxx}) holds. Furthermore,
\begin{eqnarray}\begin{aligned}[b]\label{vvvv}\|H^{p,m}_{\psi}\|_{M\dot{K}^{\alpha_{1},\lambda_{1}}_{p_{1},q_{1}}(\mathbb{Q}_p^n)\times...\times M\dot{K}^{\alpha_{m},\lambda_{m}}_{p_{m},q_{m}}(\mathbb{Q}_p^n)\rightarrow M\dot{K}^{\alpha,\lambda}_{p,q}(\mathbb{Q}_p^n)}=&\int_{\mathbb{{Z}}_{p}^*}...\int_{\mathbb{{Z}}_{p}^*}|\mathbf{t}_{1}|_{p}^{\alpha_{1}-\lambda_{1}-n(1-/q_{1})}...|\mathbf{t}_{m}|_{p}^{\alpha_{m}-\lambda_{m}-n(1-1/q_{m})}\\&\quad\times\psi(\mathbf{t}_{1},...,\mathbf{t}_{m})d\mathbf{t}_{1}...d\mathbf{t}_{m}.\end{aligned}\end{eqnarray}
\end{theorem}
The proof of Theorem (\ref{T101}) and (\ref{T11}) are more or less same. So, we just prove Theorem (\ref{T101}).

\noindent\textit{ \textbf{Proof of Theorem \ref{T101}:}}
From the previous theorem, we have:
\begin{eqnarray*}\begin{aligned}\|H^{p,2}_{\psi}(f_{1},f_{2})\chi_{k}\|_{L^{q}(\mathbb{Q}_p^n)}\leq&\int\int_{\mathbb{Z}_{p}^*}\bigg(\|f_{1}\chi_{k-m}\|_{L^{q_{1}}(\mathbb{Q}_p^n)}.\|f_{2}\chi_{k-l}\|_{L^{q_{2}}(\mathbb{Q}_p^n)}\bigg)\\\quad&\times|\mathbf{t}_{1}|_{p}^{-n/q_{1}}|\mathbf{t}_{2}|_{p}^{-n/q_{2}}\psi(\mathbf{t}_{1},\mathbf{t}_{2})d\mathbf{t}_{1}d\mathbf{t}_{2}.
\end{aligned}\end{eqnarray*}
For $1/p=1/p_{1}+1/p_{2},$ $\alpha=\alpha_{1}+\alpha_{2}$ and $\lambda=\lambda_{1}+\lambda_{2}.$ Applying H\"{o}lder's inequality together with Minkowki's inequality, we are down to:
\begin{eqnarray*}\begin{aligned}[b]\|H^{p,2}_{\psi}(f_{1},f_{2})\|_{M\dot{K}^{\alpha,\lambda}_{p,q}(\mathbb{Q}_p^n)}=&\sup_{k_{0}\in\mathbb{Z}}p^{-k_{0}\lambda}\bigg(\sum_{k=-\infty}^{k_{0}}p^{k\alpha p}\|H^{p,2}_{\psi}(f_{1},f_{2})\chi_{k}\|^{p}_{L^{q}(\mathbb{Q}_p^n)}\bigg)^{1/p}\\
\end{aligned}\end{eqnarray*}
\begin{eqnarray*}\begin{aligned}
\leq&\sup_{k_{0}\in\mathbb{Z}}p^{-k_{0}\lambda}\bigg(\sum_{k=-\infty}^{k_{0}}p^{k\alpha p}\bigg(\int\int_{\mathbb{Z}_{p}^*}\|f_{1}\chi_{k-m}\|_{L^{q_{1}}(\mathbb{Q}_p^n)}.\|f_{2}\chi_{k-l}\|_{L^{q_{2}}(\mathbb{Q}_p^n)}\\\quad&\times|\mathbf{t}_{1}|_{p}^{-n/q_{1}}|\mathbf{t}_{2}|_{p}^{-n/q_{2}}\psi(\mathbf{t}_{1},\mathbf{t}_{2})d\mathbf{t}_{1}d\mathbf{t}_{2}\bigg)^{p}\bigg)^{1/p}\\
\end{aligned}\end{eqnarray*}
\begin{eqnarray}\begin{aligned}[b]
\leq&\int\int_{\mathbb{Z}_{p}^*}\sup_{k_{0}\in\mathbb{Z}}p^{-k_{0}\lambda}\bigg(\sum_{k=-\infty}^{k_{0}}p^{k\alpha p}\bigg(\|f_{1}\chi_{k-m}\|_{L^{q_{1}}(\mathbb{Q}_p^n)}.\|f_{2}\chi_{k-l}\|_{L^{q_{2}}(\mathbb{Q}_p^n)}\bigg)^{p}\bigg)^{1/p}\\\quad&\times|\mathbf{t}_{1}|_{p}^{-n/q_{1}}|\mathbf{t}_{2}|_{p}^{-n/q_{2}}\psi(\mathbf{t}_{1},\mathbf{t}_{2})d\mathbf{t}_{1}d\mathbf{t}_{2}\\
\leq&\int\int_{\mathbb{Z}_{p}^*}\sup_{k_{0}\in\mathbb{Z}}p^{-k_{0}\lambda}\bigg(\sum_{k=-\infty}^{k_{0}}p^{k\alpha_{1} p_{1}}\|f_{1}\chi_{k-m}\|^{p_{1}}_{L^{q_{1}}(\mathbb{Q}_p^n)}\bigg)^{1/p_{1}}\\\quad&\times\bigg(\sum_{k=-\infty}^{k_{0}}p^{k\alpha_{2} p_{2}}\|f_{2}\chi_{k-l}\|^{p_{2}}_{L^{q_{2}}(\mathbb{Q}_p^n)}\bigg)^{1/p_{2}}|\mathbf{t}_{1}|_{p}^{-n/q_{1}}|\mathbf{t}_{2}|_{p}^{-n/q_{2}}\psi(\mathbf{t}_{1},\mathbf{t}_{2})d\mathbf{t}_{1}d\mathbf{t}_{2}\\
\leq&\int\int_{\mathbb{Z}_{p}^*}\sup_{k_{0}\in\mathbb{Z}}p^{-k_{0}\lambda_{1}}\bigg(\sum_{k=-\infty}^{k_{0}}p^{k\alpha_{1} p_{1}}\|f_{1}\chi_{k-m}\|^{p_{1}}_{L^{q_{1}}(\mathbb{Q}_p^n)}\bigg)^{1/p_{1}}\\\quad&\times\sup_{k_{0}\in\mathbb{Z}}p^{-k_{0}\lambda_{2}}\bigg(\sum_{k=-\infty}^{k_{0}}p^{k\alpha_{2} p_{2}}\|f_{2}\chi_{k-l}\|^{p_{2}}_{L^{q_{2}}(\mathbb{Q}_p^n)}\bigg)^{1/p_{2}}\\\quad&\times|\mathbf{t}_{1}|_{p}^{-n/q_{1}}|\mathbf{t}_{2}|_{p}^{-n/q_{2}}\psi(\mathbf{t}_{1},\mathbf{t}_{2})d\mathbf{t}_{1}d\mathbf{t}_{2}\\
\leq&\int\int_{\mathbb{Z}_{p}^*}\sup_{k_{0}\in\mathbb{Z}}p^{-(k_{0}-m)\lambda_{1}}\bigg(\sum_{k=-\infty}^{k_{0}}p^{k\alpha_{1} p_{1}}\|f_{1}\chi_{k}\|^{p_{1}}_{L^{q_{1}}(\mathbb{Q}_p^n)}\bigg)^{1/p_{1}}\\\quad&\times\sup_{k_{0}\in\mathbb{Z}}p^{-(k_{0}-l)\lambda_{2}}\bigg(\sum_{k=-\infty}^{k_{0}}p^{k\alpha_{2} p_{2}}\|f_{2}\chi_{k}\|^{p_{2}}_{L^{q_{2}}(\mathbb{Q}_p^n)}\bigg)^{1/p_{2}}\\\quad&\times|\mathbf{t}_{1}|_{p}^{-(\alpha_{1}+n/q_{1}-\lambda_{1})}|\mathbf{t}_{2}|_{p}^{-(\alpha_{2}+n/q_{2}-\lambda_{2})}\psi(\mathbf{t}_{1},\mathbf{t}_{2})d\mathbf{t}_{1}d\mathbf{t}_{2}\\
\leq&\|f_{1}\|_{M\dot{K}^{\alpha_{1},\lambda_{1}}_{p_{1},q_{1}}(\mathbb{Q}_p^n)}\|f_{2}\|_{M\dot{K}^{\alpha_{2},\lambda_{2}}_{p_{2},q_{2}}(\mathbb{Q}_p^n)}\int\int_{\mathbb{Z}_{p}^*}|\mathbf{t}_{1}|_{p}^{-(\alpha_{1}+n/q_{1}-\lambda_{1})}|\mathbf{t}_{2}|_{p}^{-(\alpha_{2}+n/q_{2}-\lambda_{2})}\\\quad&\times\psi(\mathbf{t}_{1},\mathbf{t}_{2})d\mathbf{t}_{1}d\mathbf{t}_{2}.
\end{aligned}\end{eqnarray}
First part of Theorem is done.\\
On the other hand, we define
\begin{eqnarray*}\begin{aligned}f_{1}(\mathbf{x})=|\mathbf{x}|_{p}^{-(\alpha_{1}+n/q_{1}-\lambda_{1})},\quad \mathbf{x}\in\mathbb{Q}_p^n,
\end{aligned}\end{eqnarray*}
\begin{eqnarray*}\begin{aligned}f_{2}(\mathbf{x})=|\mathbf{x}|_{p}^{-(\alpha_{2}+n/q_{2}-\lambda_{2})},\quad \mathbf{x}\in\mathbb{Q}_p^n.
\end{aligned}\end{eqnarray*}
When $\alpha_{1}\neq \lambda_{1}$ and $\alpha_{2}\neq\lambda_{2},$ we acquire
\begin{eqnarray*}\begin{aligned}\|f_{1}\chi_{k}\|_{L^{q_{1}}(\mathbb{Q}_p^n)}=&\bigg(\int_{S_{k}}|\mathbf{x}|_{p}^{-(\alpha_{1}+n/q_{1}-\lambda_{1})q_{1}}\bigg)^{1/q_{1}}\\
=&(1-p^{-n})^{1/q_{1}}p^{k(\lambda_{1}-\alpha_{1})}.
\end{aligned}\end{eqnarray*}
Also, it is not difficult to obtain
\begin{eqnarray*}\begin{aligned}\|f_{1}\|_{M\dot{K}^{\alpha_{1},\lambda_{1}}_{p_{1},q_{1}}(\mathbb{Q}_p^n)}=(1-p^{-n})^{1/q_{1}}\frac{p^{\lambda_{1}}}{(p^{p_{1}\lambda_{1}}-1)^{1/p_{1}}}.
\end{aligned}\end{eqnarray*}
Similarly, we can get
\begin{eqnarray*}\begin{aligned}\|f_{2}\|_{M\dot{K}^{\alpha_{2},\lambda_{2}}_{p_{2},q_{2}}(\mathbb{Q}_p^n)}=(1-p^{-n})^{1/q_{2}}\frac{p^{\lambda_{2}}}{(p^{p_{2}\lambda_{2}}-1)^{1/p_{2}}}.
\end{aligned}\end{eqnarray*}
For $\lambda=\lambda_{1}+\lambda_{2},$ $\alpha=\alpha_{1}+\alpha_{2}$ and $1/q=1/q_{1}+1/q_{2},$ we acquire
\begin{eqnarray*}\begin{aligned}H^{p,2}_{\psi}(f_{1},f_{2})(\mathbf{x})=|\mathbf{x}|_{p}^{-(\alpha+n/q-\lambda)}\int\int_{\mathbb{Z}_{p}^*}|\mathbf{t}_{1}|_{p}^{-(\alpha_{1}+n/q_{1}-\lambda_{1})}|\mathbf{t}_{2}|_{p}^{-(\alpha_{2}+n/q_{2}-\lambda_{2})}\psi(\mathbf{t}_{1},\mathbf{t}_{2})d\mathbf{t}_{1}d\mathbf{t}_{2}.
\end{aligned}\end{eqnarray*}
Ultimately, we have:
\begin{eqnarray*}\begin{aligned}\|H^{p,2}_{\psi}(f_{1},f_{2})\chi_{k}\|^{p}_{L^{q}(\mathbb{Q}_p^n)}=&\bigg(\int_{S_{k}}|\mathbf{x}|_{p}^{-(\alpha+n/q-\lambda)q}\bigg(\int\int_{\mathbb{Z}_{p}^*}|\mathbf{t}_{1}|_{p}^{-(\alpha_{1}+n/q_{1}-\lambda_{1})}|\mathbf{t}_{2}|_{p}^{-(\alpha_{2}+n/q_{2}-\lambda_{2})}\\\quad&\times\psi(\mathbf{t}_{1},\mathbf{t}_{2})d\mathbf{t}_{1}d\mathbf{t}_{2}\bigg)^{q}d\mathbf{x}\bigg)^{p/q}\\
=&\bigg(\int_{S_{k}}|\mathbf{x}|_{p}^{-(\alpha+n/q-\lambda)q}d\mathbf{x}\bigg)^{p/q}\bigg(\int\int_{\mathbb{Z}_{p}^*}|\mathbf{t}_{1}|_{p}^{-(\alpha_{1}+n/q_{1}-\lambda_{1})}|\mathbf{t}_{2}|_{p}^{-(\alpha_{2}+n/q_{2}-\lambda_{2})}\\\quad&\times\psi(\mathbf{t}_{1},\mathbf{t}_{2})d\mathbf{t}_{1}d\mathbf{t}_{2}\bigg)^{p}\\
=&(1-p^{-n})^{p/q}p^{-k(\alpha-\lambda)p}\bigg(\int\int_{\mathbb{Z}_{p}^*}|\mathbf{t}_{1}|_{p}^{-(\alpha_{1}+n/q_{1}-\lambda_{1})}|\mathbf{t}_{2}|_{p}^{-(\alpha_{2}+n/q_{2}-\lambda_{2})}\\\quad&\times\psi(\mathbf{t}_{1},\mathbf{t}_{2})d\mathbf{t}_{1}d\mathbf{t}_{2}\bigg)^{p}.
\end{aligned}\end{eqnarray*}
Since $\lambda_{1}=\lambda_{2}=(1/2)\lambda,$ $p_{1}=p_{2}=2p$ and $q_{1}=q_{2}=2q,$ we have:
\begin{eqnarray*}\begin{aligned}\|H^{p,2}_{\psi}(f_{1},f_{2})\|_{M\dot{K}^{\alpha,\lambda}_{p,q}(\mathbb{Q}_p^n)}=&\sup_{k_{0}\in\mathbb{Z}}p^{-k_{0}\lambda}\bigg(\sum_{k=-\infty}^{k_{0}}p^{k\alpha p}\|H^{p,2}_{\psi}(f_{1},f_{2})\chi_{k}\|^{p}_{L^{q}(\mathbb{Q}_p^n)}\bigg)^{1/p}\\
=&(1-p^{-n})^{1/q}\sup_{k_{0}\in\mathbb{Z}}p^{-k_{0}\lambda}\bigg(\sum_{k=-\infty}^{k_{0}}p^{k\alpha p}p^{-k(\alpha-\lambda)p}\\\quad&\times\bigg(\int\int_{\mathbb{Z}_{p}^*}|\mathbf{t}_{1}|_{p}^{-(\alpha_{1}+n/q_{1}-\lambda_{1})}|\mathbf{t}_{2}|_{p}^{-(\alpha_{2}+n/q_{2}-\lambda_{2})}\\\quad&\times\psi(\mathbf{t}_{1},\mathbf{t}_{2})d\mathbf{t}_{1}d\mathbf{t}_{2}\bigg)^{p}\bigg)^{1/p}\\
=&(1-p^{-n})^{1/q}\sup_{k_{0}\in\mathbb{Z}}p^{-k_{0}\lambda}\bigg(\sum_{k=-\infty}^{\infty}p^{k\lambda p}\bigg)^{1/p}\\\quad&\times\int\int_{\mathbb{Z}_{p}^*}|\mathbf{t}_{1}|_{p}^{-(\alpha_{1}+n/q_{1}-\lambda_{1})}|\mathbf{t}_{2}|_{p}^{-(\alpha_{2}+n/q_{2}-\lambda_{2})}\psi(\mathbf{t}_{1},\mathbf{t}_{2})d\mathbf{t}_{1}d\mathbf{t}_{2}\\
\end{aligned}\end{eqnarray*}
\begin{eqnarray*}\begin{aligned}
=&(1-p^{-n})^{1/q}\frac{p^{\lambda}}{(p^{\lambda p}-1)^{1/p}}\\\quad&\times\int\int_{\mathbb{Z}_{p}^*}|\mathbf{t}_{1}|_{p}^{-(\alpha_{1}+n/q_{1}-\lambda_{1})}|\mathbf{t}_{2}|_{p}^{-(\alpha_{2}+n/q_{2}-\lambda_{2})}\psi(\mathbf{t}_{1},\mathbf{t}_{2})d\mathbf{t}_{1}d\mathbf{t}_{2}\\
=&\|f_{1}\|_{M\dot{K}^{\alpha_{1},\lambda_{1}}_{p_{1},q_{1}}(\mathbb{Q}_p^n)}\|f_{2}\|_{M\dot{K}^{\alpha_{2},\lambda_{2}}_{p_{2},q_{2}}(\mathbb{Q}_p^n)}\\\quad&\times\int\int_{\mathbb{Z}_{p}^*}|\mathbf{t}_{1}|_{p}^{-(\alpha_{1}+n/q_{1}-\lambda_{1})}|\mathbf{t}_{2}|_{p}^{-(\alpha_{2}+n/q_{2}-\lambda_{2})}\psi(\mathbf{t}_{1},\mathbf{t}_{2})d\mathbf{t}_{1}d\mathbf{t}_{2}.
\end{aligned}\end{eqnarray*}
Thus (\ref{vvv}) holds in the very case.\\
When $\alpha_{1}=\lambda_{1}$ and $\alpha_{2}=\lambda_{2},$ then it is eminent to see that
\begin{eqnarray*}\begin{aligned}\|f_{1}\chi_{k}\|^{q_{1}}_{L^{q_{1}}(\mathbb{Q}_p^n)}=\|f_{2}\chi_{k}\|^{q_{2}}_{L^{q_{2}}(\mathbb{Q}_p^n)}=\int_{S_{k}}|\mathbf{x}|_{p}^{-n}d\mathbf{x}=(1-p^{-n})^{1/q_{1}}.
\end{aligned}\end{eqnarray*}
It is not hard to see that
\begin{eqnarray*}\begin{aligned}\|f_{1}\|_{M\dot{K}^{\alpha_{1},\lambda_{1}}_{p_{1},q_{1}}(\mathbb{Q}_p^n)}=p^{\lambda_{1}}(1-p^{-n})^{1/q_{1}}\frac{1}{(p^{\lambda_{1}p^{1}}-1)^{1/p_{1}}},
\end{aligned}\end{eqnarray*}
\begin{eqnarray*}\begin{aligned}\|f_{2}\|_{M\dot{K}^{\alpha_{2},\lambda_{2}}_{p_{2},q_{2}}(\mathbb{Q}_p^n)}=p^{\lambda_{2}}(1-p^{-n})^{1/q_{2}}\frac{1}{(p^{\lambda_{2}p^{2}}-1)^{1/p_{2}}}.
\end{aligned}\end{eqnarray*}
So,
\begin{eqnarray*}\begin{aligned}H^{p,2}_{\psi}(f_{1},f_{2})(\mathbf{x})=|\mathbf{x}|_{p}^{-n/q}\int\int_{\mathbb{Z}_{p}^*}|\mathbf{t}_{1}|_{p}^{-n/q_{1}}|\mathbf{t}_{2}|_{p}^{-n/q_{2}}\psi(\mathbf{t}_{1},\mathbf{t}_{2})d\mathbf{t}_{1}d\mathbf{t}_{2}.
\end{aligned}\end{eqnarray*}
Now, we have:
\begin{eqnarray*}\begin{aligned}\|H^{p,2}_{\psi}(f_{1},f_{2})\chi_{k}\|_{L^{q_{1}}(\mathbb{Q}_p^n)}=(1-p^{-n})^{1/q}\int\int_{\mathbb{Z}_{p}^*}|\mathbf{t}_{1}|_{p}^{-n/q_{1}}|\mathbf{t}_{2}|_{p}^{-n/q_{2}}\psi(\mathbf{t}_{1},\mathbf{t}_{2})d\mathbf{t}_{1}d\mathbf{t}_{2}.
\end{aligned}\end{eqnarray*}
Thus,
\begin{eqnarray*}\begin{aligned}\|H^{p,2}_{\psi}(f_{1},f_{2})\|_{M\dot{K}^{\alpha,\lambda}_{p,q}(\mathbb{Q}_p^n)}=&\sup_{k_{0}\in\mathbb{Z}}p^{-k_{0}\lambda}\bigg(\sum_{k=-\infty}^{k_{0}}p^{k\alpha p}\|H^{p,2}_{\psi}(f_{1},f_{2})\chi_{k}\|^{p}_{L^{q}(\mathbb{Q}_p^n)}\bigg)^{1/p}\\
=&(1-p^{-n})^{1/q}\sup_{k_{0}\in\mathbb{Z}}p^{-k_{0}\lambda}\bigg(\sum_{k=-\infty}^{k_{0}}p^{k\alpha p}\bigg)^{1/p}\\&\quad\times\int\int_{\mathbb{Z}_{p}^*}|\mathbf{t}_{1}|_{p}^{-n/q_{1}}|\mathbf{t}_{2}|_{p}^{-n/q_{2}}\psi(\mathbf{t}_{1},\mathbf{t}_{2})d\mathbf{t}_{1}d\mathbf{t}_{2}
\end{aligned}\end{eqnarray*}
As $\lambda_{1}=\alpha_{1}$ and $\lambda_{2}=\alpha_{2},$ we get $\lambda=\alpha.$ Hence
\begin{eqnarray*}\begin{aligned}\|H^{p,2}_{\psi}(f_{1},f_{2})\|_{M\dot{K}^{\alpha,\lambda}_{p,q}(\mathbb{Q}_p^n)}=&p^{\lambda}(1-p^{-n})^{1/q}\frac{1}{(p^{\lambda p}-1)^{1/p}}\\&\quad\times\int\int_{\mathbb{Z}_{p}^*}|\mathbf{t}_{1}|_{p}^{-n/q_{1}}|\mathbf{t}_{2}|_{p}^{-n/q_{2}}\psi(\mathbf{t}_{1},\mathbf{t}_{2})d\mathbf{t}_{1}d\mathbf{t}_{2}\\
=&\|f_{1}\|_{M\dot{K}^{\alpha_{1},\lambda_{1}}_{p_{1},q_{1}}(\mathbb{Q}_p^n)}\|f_{2}\|_{M\dot{K}^{\alpha_{2},\lambda_{2}}_{p_{2},q_{2}}(\mathbb{Q}_p^n)}\\&\quad\times\int\int_{\mathbb{Z}_{p}^*}|\mathbf{t}_{1}|_{p}^{-n/q_{1}}|\mathbf{t}_{2}|_{p}^{-n/q_{2}}\psi(\mathbf{t}_{1},\mathbf{t}_{2})d\mathbf{t}_{1}d\mathbf{t}_{2}.
\end{aligned}\end{eqnarray*}
This shows (\ref{vvv}) is also valid in this particular case.\\
Furthermore, when either $\alpha_{1}=\lambda_{1}$ or $\alpha_{2}=\lambda_{2}$ holds, we assume former holds but the later doesn't, then on the basis of previous computations we have:
\begin{eqnarray*}\begin{aligned}\|f_{1}\|_{M\dot{K}^{\alpha_{1},\lambda_{1}}_{p_{1},q_{1}}(\mathbb{Q}_p^n)}=p^{\lambda_{1}}(1-p^{-n})^{1/q_{1}}\frac{1}{(p^{\lambda_{1}p_{1}}-1)^{1/p_{1}}},
\end{aligned}\end{eqnarray*}
\begin{eqnarray*}\begin{aligned}\|f_{2}\|_{M\dot{K}^{\alpha_{2},\lambda_{2}}_{p_{2},q_{2}}(\mathbb{Q}_p^n)}=p^{\lambda_{2}}(1-p^{-n})^{1/q_{2}}\frac{1}{(p^{\lambda_{2}p_{2}}-1)^{1/p_{2}}}.
\end{aligned}\end{eqnarray*}
We definitely have the following representation
\begin{eqnarray*}\begin{aligned}H^{p,2}_{\psi}(f_{1},f_{2})(\mathbf{x})=|\mathbf{x}|_{p}^{-(\alpha_{2}+n/q-\lambda_{2})}\int\int_{\mathbb{Z}_{p}^*}|\mathbf{t}_{1}|_{p}^{-n/q_{1}}|\mathbf{t}_{2}|_{p}^{-(\alpha_{2}+n/q_{2}-\lambda_{2})}\psi(\mathbf{t}_{1},\mathbf{t}_{2})d\mathbf{t}_{1}d\mathbf{t}_{2}.
\end{aligned}\end{eqnarray*}
Now, we have:
\begin{eqnarray*}\begin{aligned}\|H^{p,2}_{\psi}(f_{1},f_{2})\chi_{k}\|_{L^{q}(\mathbb{Q}_p^n)}=&(1-p^{-n})^{1/q}p^{-k(\alpha_{2}-\lambda_{2})}\\&\quad\times\int\int_{\mathbb{Z}_{p}^*}|\mathbf{t}_{1}|_{p}^{-n/q_{1}}|\mathbf{t}_{2}|_{p}^{-(\alpha_{2}+n/q_{2}-\lambda_{2})}\psi(\mathbf{t}_{1},\mathbf{t}_{2})d\mathbf{t}_{1}d\mathbf{t}_{2}.
\end{aligned}\end{eqnarray*}
At the very end, we obtain
\begin{eqnarray*}\begin{aligned}\|H^{p,2}_{\psi}(f_{1},f_{2})\|_{M\dot{K}^{\alpha,\lambda}_{p,q}(\mathbb{Q}_p^n)}=&\sup_{k_{0}\in\mathbb{Z}}p^{-k_{0}\lambda}\bigg(\sum_{k=-\infty}^{k_{0}}p^{k\alpha p}\|H^{p,2}_{\psi}(f_{1},f_{2})\chi_{k}\|^{p}_{L^{q}(\mathbb{Q}_p^n)}\bigg)^{1/p}\\
=&(1-p^{-n})^{1/q}\sup_{k_{0}\in\mathbb{Z}}p^{-k_{0}\lambda}\bigg(\sum_{k=-\infty}^{k_{0}}p^{k\alpha p}p^{-(\alpha_{2}-\lambda_{2}kp)}\bigg)^{1/p}\\&\quad\times\int\int_{\mathbb{Z}_{p}^*}|\mathbf{t}_{1}|_{p}^{-n/q_{1}}|\mathbf{t}_{2}|_{p}^{-(\alpha_{2}+n/q_{2}-\lambda_{2})}\psi(\mathbf{t}_{1},\mathbf{t}_{2})d\mathbf{t}_{1}d\mathbf{t}_{2}\\
=&(1-p^{-n})^{1/q}p^{\lambda}\frac{1}{(p^{\lambda p}-1)^{1/p}}\\&\quad\times\int\int_{\mathbb{Z}_{p}^*}|\mathbf{t}_{1}|_{p}^{-n/q_{1}}|\mathbf{t}_{2}|_{p}^{-(\alpha_{2}+n/q_{2}-\lambda_{2})}\psi(\mathbf{t}_{1},\mathbf{t}_{2})d\mathbf{t}_{1}d\mathbf{t}_{2}.
\end{aligned}\end{eqnarray*}
Since $\alpha_{1}=\alpha_{2}=(1/2)\alpha,$ $p_{1}=p_{2}=2p,$ $q_{1}=q_{2}=2q$ and $\lambda_{1}=\lambda_{2}=(1/2)\lambda,$ we have:
\begin{eqnarray*}\begin{aligned}\|H^{p,2}_{\psi}(f_{1},f_{2})\|_{M\dot{K}^{\alpha,\lambda}_{p,q}(\mathbb{Q}_p^n)}=&\|f_{1}\|_{M\dot{K}^{\alpha_{1},\lambda_{1}}_{p_{1},q_{1}}(\mathbb{Q}_p^n)}\|f_{2}\|_{M\dot{K}^{\alpha_{2},\lambda_{2}}_{p_{2},q_{2}}(\mathbb{Q}_p^n)}\\&\quad\times\int\int_{\mathbb{Z}_{p}^*}|\mathbf{t}_{1}|_{p}^{-n/q_{1}}|\mathbf{t}_{2}|_{p}^{-(\alpha_{2}+n/q_{2}-\lambda_{2})}\psi(\mathbf{t}_{1},\mathbf{t}_{2})d\mathbf{t}_{1}d\mathbf{t}_{2}.
\end{aligned}\end{eqnarray*}
In this case (\ref{vvv}) also holds, so we conclude the proof. 

\end{document}